# EVALUATION OF THE GO METHOD WITHIN THE UPPER LIMB KINEMATICS ANALYSIS

E. ROUX[1] – S. BOUILLAND[2] –

A.-P. GODILLON-MAQUINGHEN[1] - D. BOUTTENS[2]

[1] Laboratoire d'Automatique et de Mécanique Industrielles et Humaines, U.M.R. C.N.R.S. 8530

Université de Valenciennes, Le Mont Houy, 59313 VALENCIENNES Cedex 9 - France

[2] Groupe Hopale, 62 608 BERCK SUR MER – France

**ABSTRACT**

The aim of this study is to assess the performances of the global optimisation (GO) method (Lu and O'Connor, 1999) within the upper limb kinematics analysis. First the model of the upper limb is presented. Then we apply GO method in order to reduce skin movement artefacts that imply relative movement between markers and bones. The performances of the method are then evaluated with the help of simulated movements of the upper limb. Results show a significant reduction of the errors and of the variability due to skin movement.

INTRODUCTION

Kinematics measurement techniques with external markers are commonly used within lower limb movement analysis and are more and more applied to the upper limb (Rau, et al., 2000). Markers movements relative to the underlying bones are inherent to these techniques and several methods have been proposed to reduce them. Two types of methods can be distinguished: the local or segmental methods which take into account the relative movements of the markers of a cluster attached to a body segment (Chèze, et al., 1995; Soderkvist and Wedin, 1993; Spoor and Veldpaus, 1980); the methods which optimise relative segments orientation and position thanks to joint constraints (Biryukova, et al., 2000; Lu and O'Connor, 1999; Schmidt, et al., 1999). Biryukova et al propose to optimise joint centres and axis determination but no skin movement artefacts correction is performed during voluntary movements. Schmidt et al and Lu and O'Connor compensate skin movement artefacts by controlling relative orientation and position of the segments during voluntary movements but the Global Optimisation (GO) method described by Lu and O'Connor (Lu and O'Connor, 1999) does not demand specific trials to determine the amount of skin movement artefacts to correct. GO method was initially applied to the kinematics analysis of the lower limb. We propose to assess its performances within the upper limb movement analysis thanks to simulated movements and artefacts.

METHODS

**Model**

The upper part of the human body is considered as an articulated system composed of rigid bodies corresponding to the following body segments: trunk, arm, forearm, hand (Roux, et al., 2000; Roux, et al., 2000).

*Marker locations*

Marker locations are similar to Schmidt et al's protocol (Schmidt, et al., 1998; 1999) for the acromion, the forearm, the elbow, the wrist and the hand (Fig. 1). Markers are directly attached to the subject's skin with adhesive tape. Four markers are attached to the arm. To avoid tracking difficulties due to the proximity of cluster markers, only three markers are attached to the forearm and to the hand. The trunk is characterised by markers on C7, on the 3$^{rd}$ lumbar vertebrae (L3) and on the sternum (Fig. 1).

Elbow and wrist markers are only used during a static trial since their positions, near the joints, are very sensitive to skin movements (Cappozzo, et al., 1996).

*Joint centres*

The sphere-fitting method for the determination of the joint rotation centre has proved to be more repeatable than regression methods (Leardini, et al., 1999; Stokdijk, et al., 2000). It is used to compute the rotation centre of the gleno-humeral joint, denoted shoulder centre, and the centre of the wrist. Circumduction of the shoulder and of the wrist are performed for this purpose.

The elbow centre is the middle between the medial and the lateral elbow markers.

*Definition of anatomical frames*

According to the static trial, anatomical frames are defined in regards to the ISB recommendations (Wu and Cavanagh, 1995) and are described in Table 1. Euler's angles have been chosen to describe the relative movement of the body segments.

*Joint constraints*

The elbow and the wrist joints are considered as cardanic joints. The abduction-adduction of the elbow and the pronation-supination of the wrist are forced to be within the interval [-1°, 1°]. Dislocation of the multi-link system is prevented by imposing a 2 mm maximal translation between the arm and the forearm and between the forearm and the hand. The [-1°, 1°] interval and the 2 mm translation define joints laxity. No constraint is imposed on the location of the shoulder centre. Indeed the scapular motion can not be reliability determined with the help of external markers and a motion capture system (Pronk, 1991).

**Global Optimisation applied to the upper limb**

In order to minimise relative movement between clusters and bones, we apply the GO method presented by Lu and O'Connor (Lu and O'Connor, 1999). Markers positions at the static trial are free from skin movement artefacts and are taken as reference (Lu and O'Connor, 1999). Moreover, during movement, we consider as negligible the whole displacement of the cluster of the hand relative to the underlying bones.

The weighting matrix *W* described by Lu and O'Connor is defined with segmental residual errors given by the algorithm of Söderkvist (Soderkvist and Wedin, 1993). An iterative optimisation method is then used to compute the optimal parameters.

**Evaluation of the method**

Evaluation of the GO method was carried out with simulated movements.

The static trial of a subject gave the geometrical model. Two movements were simulated: a pure internal-external rotation of the shoulder and a pure pronation-supination of the elbow. Indeed the relative movement between clusters and underlying bones especially affects the evaluation of the axial rotation (Cappozzo, et al., 1996; Roux, et al., 2000; Schmidt, et al., 1999).

*External-internal rotation of the shoulder*

The movement took the following form:

$Rot = Rot_{ArmInit} + \frac{\pi}{3} sin(2\pi\ 0.5\ t)$, with $Rot_{ArmInit}$ the initial rotation angle, i.e. the angular configuration during the static trial.

The amplitude $\frac{\pi}{3}$ ensured a realistic motion range and the frequency of 0.5 Hz a feasible movement velocity. Movement duration was two seconds, corresponding to one period of the sinusoidal wave, with a 50 Hz sample frequency.

Concerning measurement errors, we imposed a maximal error of 5.5 mm for relative distance between markers. This value is based on the results of Richards (Richards, 1999). To this end, measurement errors were considered as a random noise with a normal distribution (mean = 0 mm; standard deviation = 0.615 mm). This distribution ensured a maximum position measurement error of 1.59 mm (99% confidence) for each direction and consequently a maximum error of 5.5 mm for the distance between two markers.

The error previously described was assumed to be the worst we could find in our application and was applied to each marker position.

Skin movement artefacts were simulated by a continuous noise model of the form $A.\sin(\omega.t+\varphi)$ (Chèze, et al., 1995; Chèze, et al., 1998; Lu and O'Connor, 1999). $A$ is the amplitude of the noise, $\omega$ its frequency and $\varphi$ its phase angle. This noise was only applied to the markers of the moving segment, i.e. the arm. $A$ was assumed to be proportional to movement amplitude. This assumption was verified by Schmidt et al during a pure axial rotation of the forearm (Schmidt, et al., 1998; 1999). Given the fact that skin movements are greater near the proximal end of a segment and that displacements of the markers with respect to the underlying bone can reach 40 mm on the lower limb (Cappozzo, et al., 1996), $A$ was scaled to be between 0 and 20 mm for the two proximal markers of the arm, and between 0 and 10 mm for the two others. So for the marker $m$, the amplitude of the noise was of the form: $A_m = B_m.|\sin(2\pi\,0.5\,t)|$, with $B_m \in \{10, 20\}$.

$\omega_m$ and $\varphi_m$ are random scaled numbers. $\omega_m$ was scaled to be between $\pi$ and $3\pi$, i.e. between 1 and 3 times the frequency of the movement. The lower limit ensured that the maximum signal to noise ratio was observed during the trial. $\varphi_m$ was scaled to be between 0 and $2\pi$ (Chèze, et al., 1995; Chèze, et al., 1998; Lu and O'Connor, 1999).

*Pronation-supination of the elbow*

The same method was applied to simulate a pure pro-supination of the forearm. The simulated movement was of the form $Rot_{Forearm} = Rot_{Forearm\ Init} + \frac{\pi}{3}\sin(2\pi\,0.5\,t)$, with $Rot_{Forearm\ Init}$ the initial pro-supination angle that was imposed to be neutral.

Skin movement artefacts were added on the three markers of the forearm. $B_m = 20$ mm for the proximal marker and $B_m = 10$ mm for the two others. In practice, a pro-supination movement is combined with a contraction of the biceps muscle, implying artefacts on the position of the arm markers. So simulated skin movement artefacts were added to these markers too, with a maximum amplitude $B_m = 5$ mm.

RESULTS

Thirty simulated trials were performed. They differ on the noise parameters $\omega_m$ and $\varphi_m$, that are randomly settled in a given interval previously defined. For each variable and each trial, we computed the Root Mean Square (RMS) of the errors with and without using GO. The average and the standard deviation of the RMS over the thirty trials were computed. Figure 2 and 3 correspond to the angle errors and Figures 4 and 5 to the relative translation between adjacent body segments, denoted dislocation.

DISCUSSION

Only the relative movement between the hand and the forearm during the internal-external rotation of the shoulder present greater errors with the application of the GO method than without it. Moreover these errors still remain very low in regards to the measurement ones (Richards, 1999). Errors are significantly compensated for all the other degrees of freedom.

The variability of the results are also significantly reduced with the GO method, except for the flexion-extension and the abduction-adduction of the wrist during internal-external rotation of the shoulder.

In spite of the good results thanks to the application of the GO method, some aspects of the kinematics model are arguable. From an anatomical point of view, the flexion axis of the elbow is not normal to the plane formed by the longitudinal axis of the arm and of the forearm during the static trial. However according to Wang et al. (Wang, et al., 1998) this definition can be used to approximate the elbow flexion extension axis in a large motion range of the elbow joint, especially when the forearm is pronated and neutral. By using the joint co-ordinate system of the elbow, Schmidt et al. (Schmidt, et al., 1998; 1999) make the same assumption for the whole movement. The kinematics constraint that imposes the abduction-adduction of the elbow to be within the interval [-1°, 1°] is not verified from an anatomical point of view because of the bony structure of the forearm and of the articular surfaces of the elbow joint. However, simulations show that an abnormal abduction-adduction of the elbow is observed because of the skin movement artefacts, when GO is not applied (Figures 2 and 4). Consequently the abduction-adduction of the elbow observed during a voluntary movement and without applying the GO method is due to the skin movement artefacts and does not reflect the real bone movements.

## Conclusion

This study shows that GO method significantly reduces the errors and the variability introduced by skin movements within the kinematics analysis of the upper limb with external markers. An other advantage of the method is that it does not demand specific trials realisation to previously estimate skin movements artefacts.

**Acknowledgments**

The authors thank Olivier REMY-NERIS, MD, PhD, responsible for the *Laboratoire d'Analyse du Mouvement du Groupe Hopale,* for his assistance in performing this study.

This work is sponsored by the *Centre National de Recherche Scientifique (C.N.R.S.)*, the *Nord-Pas de Calais* Region and the *Institut Régional de Recherche sur le Handicap (I.R.R.H.).*

# References

Biryukova, E.V., Roby-Brami, A., Frolov, A.A., Mokhtari, M., 2000. Kinematics of human arm reconstructed from spatial tracking system recordings. Journal of Biomechanics 33, 985-995

Cappozzo, A., Catani, F., Leardini, A., Benedetti, M.G., Della Croce, U., 1996. Position and orientation in space of bones during movement: experimental artefacts. Clinical Biomechanics 11, 90-100

Chèze, L., Fregly, B.J., Dimnet, J., 1995. A solidification procedure to facilitate kinematics analyses based on video system data. Journal of Biomechanics 28, 879-884

Chèze, L., Fregly, B.J., Dimnet, J., 1998. Determination of joint functional axes from noisy marker data using the finite helical axis. Human Movement Science 17, 1-15

Leardini, A., Cappozzo, A., Catani, F., Toksvig-Larsen, S., Petitto, A., Sforza, V., Cassanelli, G., Giannini, S., 1999. Validation of a functional method for the estimation of hip joint centre location. Journal of Biomechanics 32, 99-103

Lu, T.-W., O'Connor, J.J., 1999. Bone position estimation from skin marker co-ordinates using global optimisation with joint constraints. Journal of Biomechanics 32, 129-134

Pronk, G.M. (1991) The shoulder girdle. PhD. thesis, Delft Univ. Press, Delft

Rau, G., Disselhorst-Klug, C., Schmidt, R., 2000. Movement biomechanics goes upwards: from the leg to the arm. J Biomech 33, 1207-1216.

Richards, J.G., 1999. The measurement of human motion: A comparison of commercially available systems. Human Movement Science 18, 589-602

Roux, E., Bouilland, S., Bouttens, D., Istas, D., Godillon-Maquinghen, A.-P., Lepoutre, F.-X., 2000. Evaluation of the kinematics of the shoulder and of the upper limb. In Proceedings of the 3rd conference of the International Shoulder Group. Newcastle Upon Tyne, U.K.

Roux, E., Bouilland, S., Bouttens, D., Istas, D., Godillon-Maquinghen, A.-P., Lepoutre, F.-X., 2000. Experimental protocol for the kinematics measurement of the shoulder and of the upper limb. In Proceedings of the Proceedings / Actes de conférence - XXVe Congrès de la Société de Biomécanique - XIe Congrès de la Société Canadienne de Biomécanique - Archives of Physiology and Biochemistry.

Schmidt, R., Disselhorst-Klug, C., Silny, J., Rau, G., 1998. A measurement procedure for the quantitative analysis of the free upper-extremity movements. In Proceedings of the Fifth international symposium on the 3-D analysis of human movement. Chattanooga, Tennessee, USA

Schmidt, R., Disselhorst-Klug, C., Silny, J., Rau, G., 1999. A marker-based measurement procedure for unconstrained wrist and elbow motions. Journal of Biomechanics 32, 615-621

Soderkvist, I., Wedin, P.A., 1993. Determining the movements of the skeleton using well-configured markers. Journal of Biomechanics 26, 1473-1477

Spoor, C.W., Veldpaus, F.E., 1980. Rigid body motion calculated from spatial co-ordinates of markers. Journal of Biomechanics 13, 391-393

Stokdijk, M., Nagels, J., Rozing, P.M., 2000. The glenohumeral joint rotation centre in vivo. Journal of Biomechanics 33, 1629-1636.

Wang, X., Maurin, M., Mazet, F., De Castro Maia, N., Voinot, K., Jean-Pierre, V., Fayet, M., 1998. Three-dimensional modelling of the motion range of axial rotation of the upper arm. Journal of Biomechanics 31, 899-908

Wu, G., Cavanagh, P.R., 1995. ISB recommendations for standardization in the reporting of kinematic data. J Biomech 28, 1257-1261.

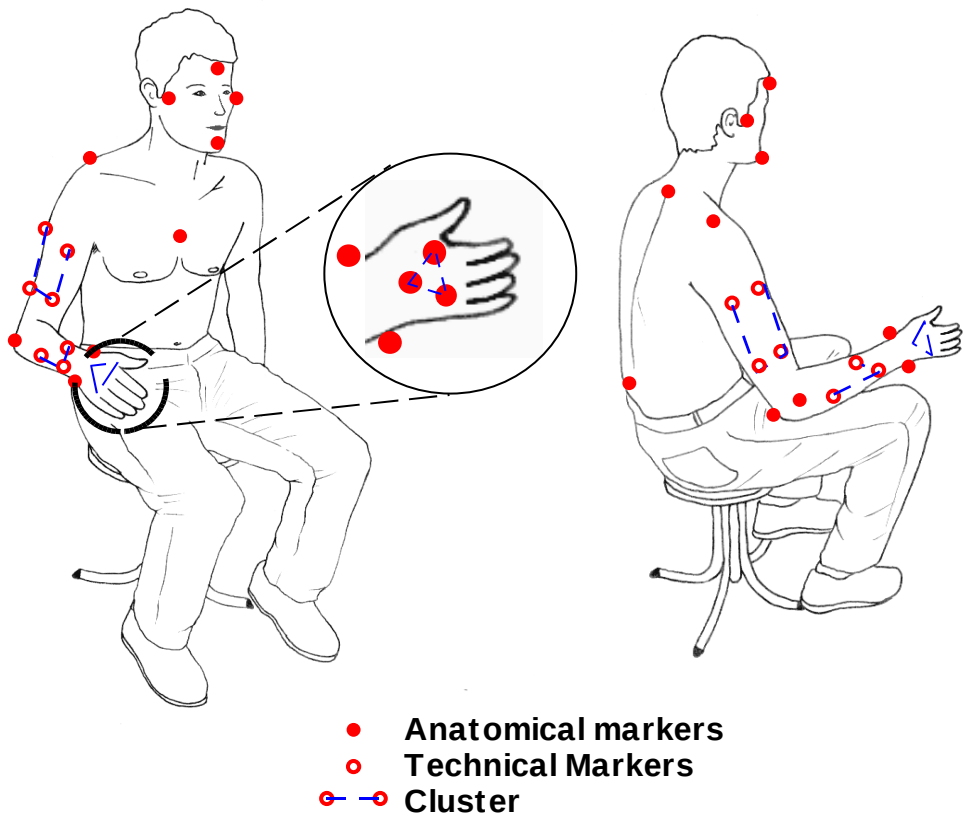

*Figure 1: Marker locations during static trial*

| SEGMENT | X-Axis | Y-Axis | Z-Axis |
|---|---|---|---|
| **Head** | Y ∧ Z | chin --> forehead | Left --> right temples |
| **Trunk** | Y ∧ Z | L3 --> C7 | (L3 --> sternum) ∧ Y |
| **Shoulder girdle** | (L3 --> C7) ∧ Z | Z ∧ X | C7 --> acromion |
| **Arm** | Y ∧ Z | Elbow centre --> shoulder centre | (longitudinal axis of the forearm) ∧ (Y-axis of the forearm) |
| **Forearm** | Y ∧ Z | Wrist centre --> elbow centre | (Posterior --> anterior styloïd) ∧ Y |
| **Hand** | Y ∧ Z | Barycentre of the three markers of the hand --> marker near the wrist | (Posterior --> anterior hand markers) ∧ Y |

*Table 1: Definition of the anatomical frames for each body segment*

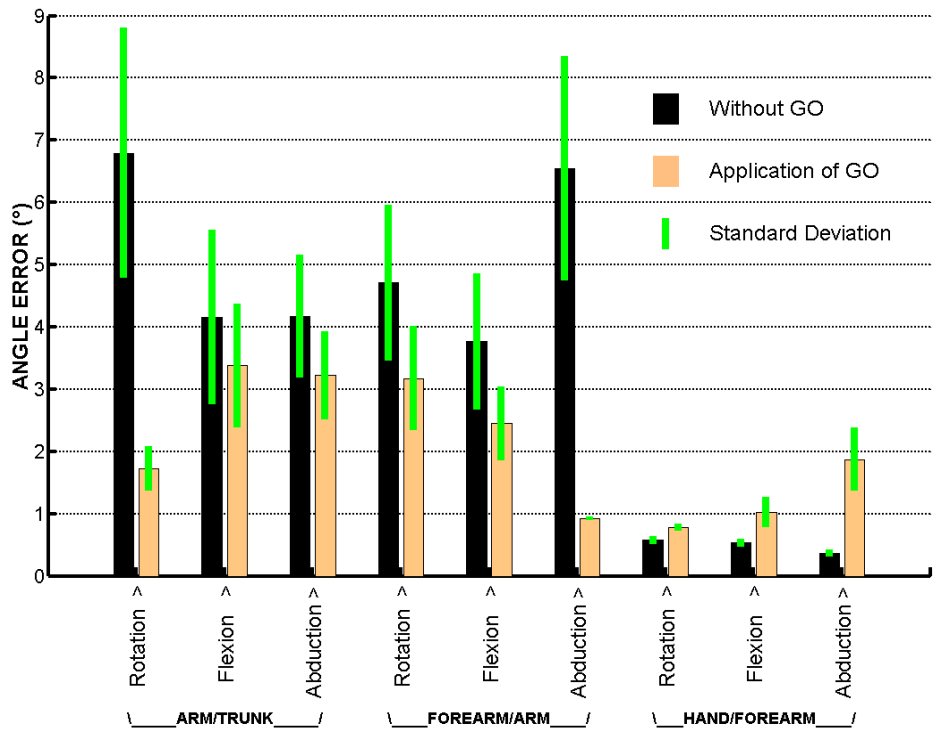

*Figure 2: Mean / Standard deviation of RMS error of angles for external-internal rotation*

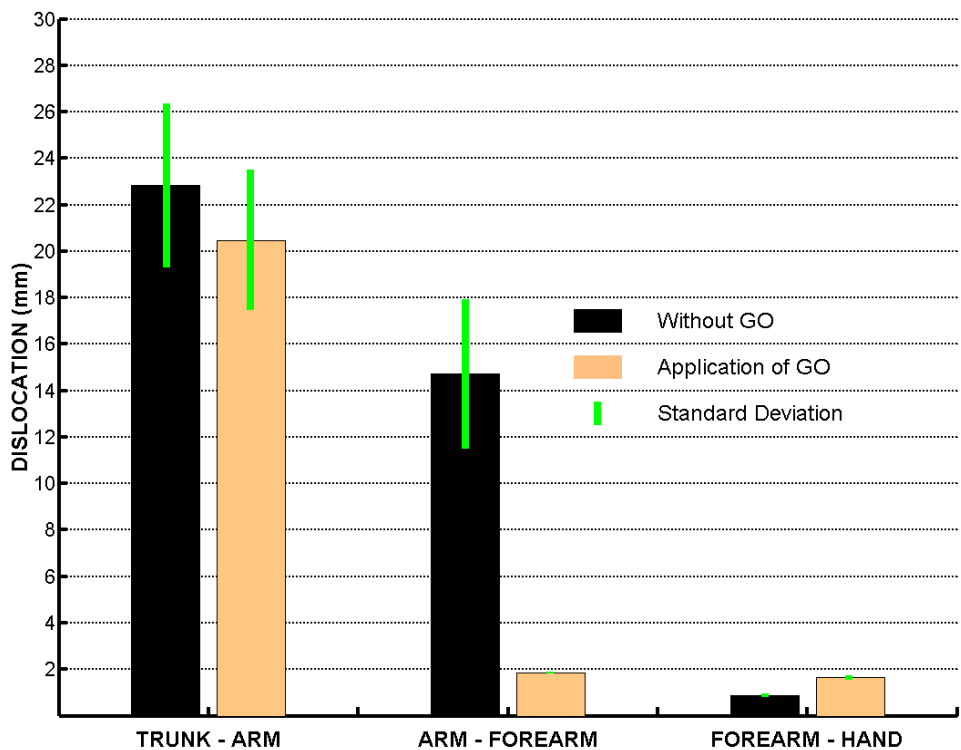

*Figure 3: Mean / Standard deviation of RMS of translation between segments for external-internal rotation*

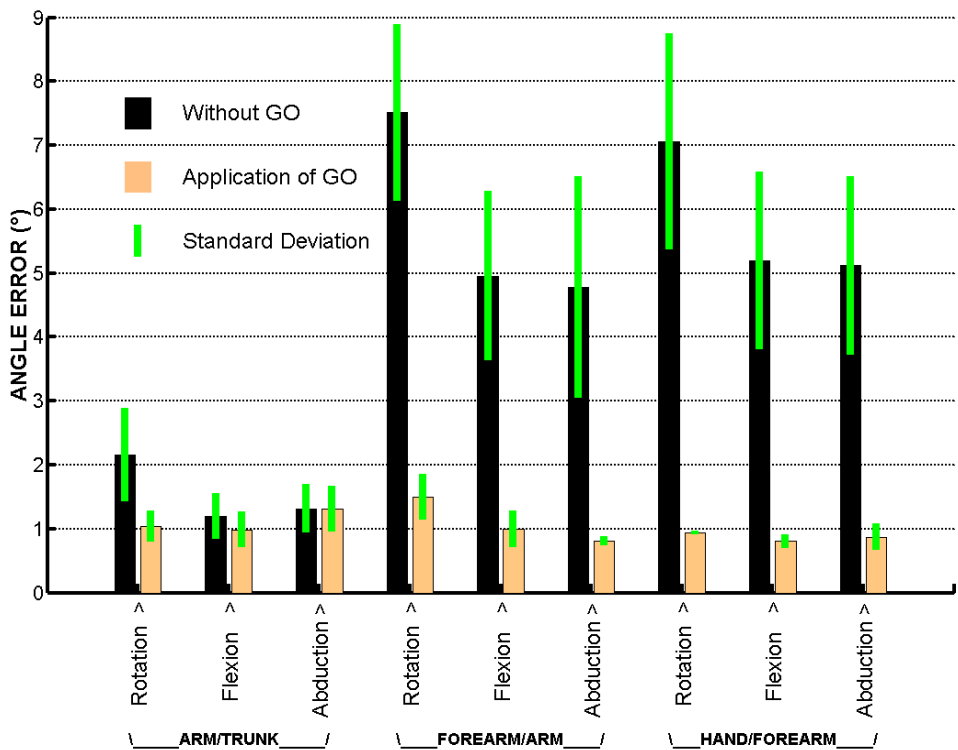

*Figure 4: Mean / Standard deviation of RMS error of angles for pro-supination of the elbow*

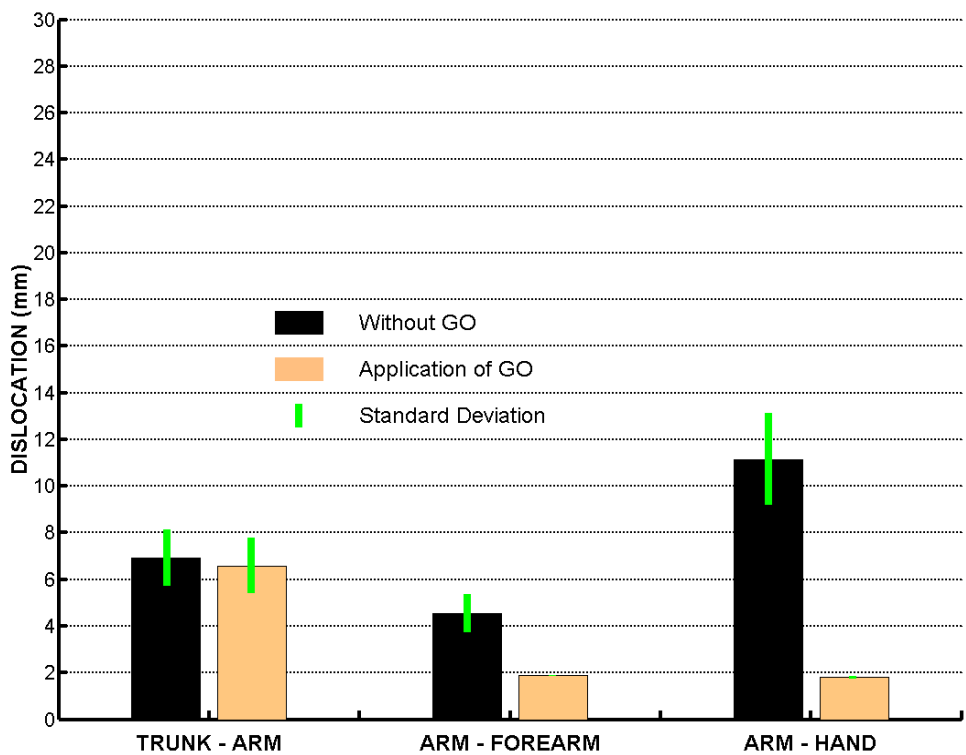

*Figure 5: Mean / Standard deviation of RMS of translation between segments*

*for pro-supination of the elbow*